\documentclass[12pt]{amsart}

\usepackage{color}
\usepackage{amstext}
\usepackage{amsfonts}
\usepackage{amsthm}
\usepackage{amsmath}
\usepackage{amssymb}
\usepackage{latexsym}
\usepackage{amsfonts}
\usepackage{graphicx}
\usepackage{texdraw}
\usepackage{graphpap}
\usepackage{subfigure}
\usepackage{enumerate}
\usepackage[pagebackref,hypertexnames=false, colorlinks, citecolor=red, linkcolor=red]{hyperref} 
\usepackage[backrefs]{amsrefs}
\usepackage{amsmath,amstext,amsthm,amssymb,amsxtra}
\usepackage{mathtools}
\DeclarePairedDelimiter{\ceil}{\lceil}{\rceil}

\usepackage{soul}
\usepackage[margin=15mm]{geometry}
\usepackage{tikz}
\usetikzlibrary{3d,calc,fadings,decorations.pathreplacing}

\input txdtools

\begin{document}

\makeatletter
\newcommand*{\defeq}{\mathrel{\rlap{%
                     \raisebox{0.3ex}{$\m@th\cdot$}}%
                     \raisebox{-0.3ex}{$\m@th\cdot$}}%
                     =}
										
\newcommand*{\eqdef}{=
										 \mathrel{\rlap{%
                     \raisebox{0.3ex}{$\m@th\cdot$}}%
                     \raisebox{-0.3ex}{$\m@th\cdot$}}%
										}
\makeatother

\newcommand{\unit}{1\!\!1}

\newcommand*\rfrac[2]{{}^{#1}\!/_{#2}}

\newcommand{\norm}[1]{\ensuremath{\left\|#1\right\|}}
\newcommand{\abs}[1]{\ensuremath{\left\vert#1\right\vert}}
\newcommand{\ip}[2]{\ensuremath{\left\langle#1,#2\right\rangle}}
\newcommand{\sgn}{\text{sgn}}

\newcommand{\pbar}{\ensuremath{\bar{\partial}}}
\newcommand{\db}{\overline\partial}
\newcommand{\D}{\mathbb{D}}
\newcommand{\B}{\mathbb{B}}
\newcommand{\Sp}{\mathbb{S}}
\newcommand{\T}{\mathbb{T}}
\newcommand{\R}{\mathbb{R}}
\newcommand{\Z}{\mathbb{Z}}
\newcommand{\C}{\mathbb{C}}
\newcommand{\N}{\mathbb{N}}

\newcommand{\m}[1]{\mathcal{#1}}
\newcommand{\eps}{\epsilon}
\newcommand{\he}[1]{h_{#1}^{\eps}}
\newcommand{\hn}[1]{h_{#1}^{\eta}}
\newcommand{\avg}[1]{\langle #1 \rangle}

\newcommand{\mcd}{\mathcal{D}}
\newcommand{\mcr}{\mathcal{R}}
\newcommand{\mch}{\mathcal{H}}
\newcommand{\mcc}{\mathcal{C}}
\newcommand{\prb}{\mathbb{P}}
\newcommand{\Prb}{\mathbb{P}}

\newcommand{\ep}{\epsilon}
\newcommand{\lb}{\lambda}

\newcommand{\La}{\langle }
\newcommand{\Ra}{\rangle }

\newcommand{\bmo}{\textnormal{BMO}}
\newcommand{\wbmod}{BMO_{\mcd}(w)}
\newcommand{\wbmosd}{BMO^2_{\mcd}(w)}
\newcommand{\nbmod}{BMO_{\mcd}(\nu)}
\newcommand{\nbmosd}{BMO^2_{\mcd}(\nu)}
\newcommand{\avgb}{\left<b\right>}
\newcommand{\h}[1]{h_{#1}^{0}}
\newcommand{\nch}[1]{h_{#1}^{1}}
\newcommand{\fd}{I_{\alpha}^{\m{D}}}
\newcommand{\BMO}{\textnormal{BMO}}
\newcommand{\pcwise}[5]{#1 = \left\{\begin{array}{lr} #2, &\hspace{1mm} #3\\ #4, &\hspace{1mm} #5 \end{array} \right.}

\numberwithin{equation}{section}

\newtheorem{thm}{Theorem}[section]
\newtheorem{lm}[thm]{Lemma}
\newtheorem{cor}[thm]{Corollary}
\newtheorem{conj}[thm]{Conjecture}
\newtheorem{prob}[thm]{Problem}
\newtheorem{prop}[thm]{Proposition}
\newtheorem*{prop*}{Proposition}
\newtheorem{defn}[thm]{Definition}
\newtheorem{remar}[thm]{Remark}

\theoremstyle{remark}
\newtheorem{rem}[thm]{\textbf{Remark}}
\newtheorem*{rem*}{\textbf{Remark}}

\title[1-Bit RIP Embeddings]
{Noisy 1-Bit Compressed Sensing Embeddings \\ Enjoy a Restricted Isometry Property}

%
%
\thanks{Research supported in part by grant NSF-DMS 1265570. }

\author{Scott Spencer}
\address{ School of Mathematics, Georgia Institute of Technology, Atlanta GA 30332, USA}
\email{spencer@math.gatech.edu}


\begin{abstract} 
We investigate the sign-linear embeddings of 1-bit compressed sensing given by Gaussian measurements.  One can give short arguments concerning a Restricted Isometry Property of such maps using \textit{Vapnik-Chervonenkis dimension} of sparse hemispheres.  This approach has a natural extension to the presence of additive white noise prior to quantization.  Noisy one-bit mappings are shown to satisfy an RIP when the metric on the sphere is given by the noise.  

\end{abstract}

\maketitle
\setcounter{tocdepth}{1}

\section{Introduction}
Compressed sensing is a modern data processing scheme that is proving useful in many scientific areas, such as MR imaging, radar, astronomy: see \cite{radar, astro, mri} for more details.  The overarching goal is to reconstruct a \textit{signal} $x\in \mathbb{R}^n$ from the \textit{measurements} $Ax\in \mathbb{R}^m$ ($m \ll n$) given the sensing matrix $A\in \mathbb{R}^{m\times n}$ and some constraint on the set of signals.  Without such a constraint, this is an ill-posed inverse problem, while more information about the signal $x$ may make the objective approachable.  One common situation is that the signal is sparse:  for a signal $x=(x^1,\dots,x^n)$, we say $x$ is $s$-sparse if $\left|\{x^j\not= 0\}\right|\leq s$.  A successful program for reconstructing sparse signals is $\ell_1$-minimization.  This convex optimization algorithm is tractable and perfectly reconstructs $s$-sparse vectors (and well approximates them in the presence of noise) if the sensing matrix $A$ has the $(s,\delta)$-RIP with small enough $\delta$ \cite{Candes}.  A matrix $A$ is said to have the $(s,\delta)$-RIP if \[(1-\delta)\|x-y\|_2^2\leq \|Ax-Ay\|_2^2\leq (1+\delta)\|x-y\|_2^2\] for all pairs $x,y$ of $s$-sparse vectors.  The object of our interest is the analogue, i.e., dimension reducing quasi-isometric embeddings of sparse vectors, in the 1-bit sensing framework.

\subsection{1-Bit Sensing}
We study the dimension reducing \textit{sign-linear} maps of 1-bit compressed sensing.  Associated to each $A\in \R^{m\times n}$ is the sign-linear map 
\begin{align}\label{E:phia}
\Phi_A:\mathbb{S}^{n-1} \to \mathcal{H}^m \\
\nonumber \Phi_Ax=\sgn(Ax),
\end{align}
where $\mathcal{H}^m$ is the Hamming Cube $\{\pm 1\}^m$, the \textit{sgn} map is applied component-wise, and \[\pcwise{\sgn(x)}{+1}{x>0}{-1}{x\leq 0.}\]  We restrict our attention to the sphere since any two signals that differ only in norm will have identical measurements.  In the larger realm of compressed sensing, 1-bit sensing is the case of extreme quantization: only the sign-bit of each linear measurement is preserved.  The concept was initially suggested by Boufounos-Baraniuk \cite{boubar} in 2008.

Let $\Sp^{n-1}_s$ denote the set of $n$-dimensional, unit length $s$-sparse signals.  The \textit{$(s,\delta)$-Restricted Isometry Property}, or $(s,\delta)$-RIP, analogue for $\Phi_A$ that we investigate is \[\sup_{x,y\in \Sp_s^{n-1}} \left|d_{\mathcal{H}^m}(\Phi_Ax,\Phi_Ay)-d(x,y)\right|\leq \delta,\] where $d(\cdot,\cdot)$ is geodesic distance on the sphere, and $d_{\mathcal{H}^m}(\cdot,\cdot)$ is the Hamming metric: \[d_{\mathcal{H}^m}(a,b):=\tfrac{1}{m}\left|\{1\leq k \leq m:a_k \not= b_k\}\right|.\]

The reader may notice that the 1-bit RIP given above is single-scale, while the original RIP is multiscale.  This modification is unavoidable; given $A\in \R^{m\times n}$ and $\eps>0$, there are $x,y \in \Sp^{n-1}_s$ such that $d(x,y)\leq \eps$ and $d_{\mathcal{H}^m}(\Phi_Ax,\Phi_Ay)\geq \frac{1}{m}.$  This formulation of the RIP has been studied theoretically, see \cite{BL,PV}; it also plays a role in sparse signal recovery from 1-bit measurements, e.g. \cite{JacquesD, binarystable}.

It is natural to consider the effects of noise on a 1-Bit embedding.  We consider the case of additive white noise prior to quantization.  When we consider a random sensing matrix $A$ and random noise $\eta$, we always assume they are independent.  Associated to a matrix $A\in\R^{m\times n}$ and vector $\eta\in\mathbb{R}^m$ is a 1-bit embedding of the form
\begin{align}\label{E:phiaeta}
\Phi^\eta_A &:\mathbb{R}^n \to \mathcal{H}^m \\
\nonumber \Phi^\eta_A x &=\sgn(Ax+\eta).
\end{align}
Taking $\eta \sim\mathcal{N}(\textbf{0},\sigma^2 I_m),$ the $(s,\delta)$-RIP analogue for $\Phi^\eta_A$ that we investigate is \[\sup_{x,y\in \Sp_s^{n-1}} \left|d_{\mathcal{H}^m}(\Phi^\eta_Ax,\Phi^\eta_Ay)-d^\sigma(x,y)\right|\leq \delta.\]  The distance $d^\sigma(\cdot,\cdot)$ is a distorted version of the geodesic distance, defined in \eqref{e:xdsigma} and discussed in Section~\ref{S:noisy}.  The affects of the additive white noise on the RIP are analyzed by increasing the Gaussian measurements' dimension by one and lifting the sphere to one higher dimension by padding with $\sigma^2$.

\section{Outline and Main Results}\label{S:mainresults}
For $x\in \Sp^{n-1}$, set $H_x =\{p\in\Sp^{n-1}:\avg{p,x}>0\}$, the hemisphere associated to $x$.  Denote by $H^{n,s}$ the set of hemispheres of $\Sp^{n-1}$ associated to $s$-sparse signals: $H^{n,s}=\{H_x:x\in\Sp_s^{n-1}\}.$  The first result listed here gives a useful upper bound on the $VC$-dimension, defined in Section~\ref{S:VCdef}, of $H^{n,s}$.  The result easily applies to half-spaces, a well studied classification scheme in learning theory; it is well known that the $VC$-dimension of half-spaces in $\R^n$ indexed by $s$-sparse vectors is $\mathcal{O}(s\log n)$.  The theorem below is slightly better, but we are unsure if it is known.  We include the proof in Section~\ref{S:VC} for completeness, and note that it is quite surprising to find the popular $s\log(n/s)$ quantity.  Throughout, $x\lesssim y$ means there is an absolute $C>0$ such that $x\leq Cy$.
\begin{thm}\label{T:mainthmvc}
$VC(H^{n,s})\lesssim s\log(n/s).$
\end{thm}

\begin{defn}\label{T:def1}
Let $\Phi:\Sp_s^{n-1}\to\mathcal{H}^m$.  We say $\Phi$ has the $(s,\delta)$-RIP if
\[\sup\limits_{x,y\in \mathbb{S}^{n-1}_s}\left| d_{\mathcal{H}^m}(\Phi x,\Phi y)-d(x,y)\right| \leq \delta.\]
\end{defn}

Of note in Definition~\ref{T:def1} is the metric $d(\cdot,\cdot)$, which is \textbf{not} the euclidean distance, but rather the geodesic distance on the sphere, normalized so that antipodal points are unit distance apart:  \[d(x,y):=\tfrac{1}{\pi}\arccos\big(\avg{x,y}\big).\]  This choice of metric is natural since it is the expectation of $d_{\mathcal{H}^m}(\Phi_Ax,\Phi_Ay)$.

In Section~\ref{S:noiseless} we employ a standard entropy integral argument to bound a supremum, indexed by pairs of $s$-sparse vectors.  This is an alternative proof of a recent result of Bilyk-Lacey, the case of sparse vectors in \cite[Theorem 1.14]{BL}, which is:


\begin{thm}\label{T:mainnoiseless2}
Let $A\in \mathbb{R}^{m\times n}$ with rows drawn independently from the standard Gaussian distribution.  Then for any $0<\varepsilon, \delta<1$ and $1\leq s<n$, $\Phi_A$ has the $(s,\delta)$-RIP with probability at least $1-\varepsilon$ provided \[m\gtrsim \delta^{-2}\left[\log(2/\varepsilon)+s\log(n/s)\right].\]
\end{thm}

The next theorem, proved in Section~\ref{S:noisy}, is the import of the paper.  We consider the 1-bit sign-linear maps with additive white noise prior to quantization.  A curious detail about the result is that the error due to noise is not naturally expressed in the distortion parameter, nor the number of measurements or probability of success, but rather in the metric on the sphere.  That is, if the sphere is endowed with a certain ``distorted'' geodesic metric \eqref{e:xdsigma}, the noisy embedding has the $(s,\delta)$-RIP with the same order of measurements and probability of success as determined in Theorem~\ref{T:mainnoiseless2}.  Before stating the theorem, we define:

\begin{equation}\label{e:xdsigma}
d^\sigma(x,y):=\tfrac{1}{\pi}\arccos\left(\tfrac{\avg{x,y}+\sigma^2}{1+\sigma^2}\right) \,\,\, \text{for}\,\, x,y\in\Sp^{n-1}.
\end{equation}

We also define the following noisy version of the 1-bit RIP:

\begin{defn}
Let $\Phi:\Sp_s^{n-1}\to\mathcal{H}^m$.  We say $\Phi$ has the $(s,\delta,\sigma)$-RIP if
\[\sup\limits_{x,y\in \mathbb{S}^{n-1}_s}\left| d_{\mathcal{H}^m}(\Phi x,\Phi y)-d^\sigma(x,y)\right| \leq \delta.\]
\end{defn}

\begin{thm}\label{T:mainnoise}
Let $A\in \mathbb{R}^{m\times n}$ with rows drawn independently from the standard Gaussian distribution and $\eta \sim \mathcal{N}(\textbf{0},\sigma^2 I_m)$.  Then for any $0<\varepsilon, \delta<1$ and $1\leq s<n$, $\Phi^\eta_A$ has the $(s,\delta,\sigma)$-RIP with probability at least $1-\varepsilon$ provided \[m\gtrsim \delta^{-2}\left[\log(2/\varepsilon)+s\log(n/s)\right].\]
\end{thm}

\begin{rem*}
It is a common goal in noisy compressive sensing to ``eliminate'' the noise.  That is, one wishes to take enough measurements so that the noise is practically negligible.  Theorem~\ref{T:mainnoise} demonstrates that this possibility is controlled by the variance in the Gaussian noise model.  The empirical process of interest approaches the distorted metric $d^\sigma$, which is a deterministic object that necessarily deviates from the geodesic metric when $\sigma^2>0$.  
\end{rem*}

We conclude with Section~\ref{S:quieting}, comparing the geodesic distance with the metric defined in \eqref{e:xdsigma}.  A crude upper bound on their difference gives a lower bound on the number of Gaussian measurements needed for a noisy embedding to have the RIP into the Hamming cube with the \emph{geodesic} metric prescribed to the sphere.  While this result is appealing for obvious reasons, Theorem~\ref{T:mainnoise} may be more useful in practice, allowing the reader to appeal to the fact that the two metrics are indeed \textit{very} close at small scales.

\begin{cor}\label{T:quieting}
Let $A\in \mathbb{R}^{m\times n}$ with rows drawn independently from the standard Gaussian distribution and $\eta \sim \mathcal{N}(\textbf{0},\sigma^2 I_m)$.  Then for any $0<\varepsilon<1$ and  $\delta>1-\frac{1}{\pi}\arccos\left(\frac{\sigma^2-1}{\sigma^2+1}\right),$ $\Phi_A^\eta$ has the $(s,\delta)$-RIP with probability at least $1-\varepsilon$ provided \[m\gtrsim \left[\delta +\tfrac{1}{\pi}\arccos\left(\tfrac{\sigma^2-1}{\sigma^2+1}\right)-1\right]^{-2}\left[\log(2/\varepsilon)+s\log(n/s)\right].\]
\end{cor}


\section{The VC-Dimension of Sparse Hemispheres}\label{S:VCsection}

\subsection{VC Dimension}\label{S:VCdef}

Let $X$ be a set and $\mathcal{C}$ be a collection of subsets of $X$.  Denote by $\binom{X}{k}$ the set of subsets of $X$ with $k$ elements.  For each $k\in\N$, define \[m^\mcc(k):= \max\limits_{B\in {X \choose k}} \left|\{B\cap C:C\in\mathcal{C}\}\right|.\]  Clearly $m^\mcc(k)\leq 2^k$.  The \textit{Vapnik-Chervonenkis dimension} ($VC-$dimension) of $\mcc$, denoted $VC(\mathcal{C})$, is the largest integer $d$ (if it exists) such that $m^\mcc(d)=2^d$, and $VC(\mathcal{C})=\infty$ otherwise.  Alternatively, we say $\mcc$ \textit{shatters} $B$ if every subset of $B$ is realized as the intersection of $B$ with an element of $\mcc$.  Then $VC(\mcc)$ is the cardinality of the largest subset it shatters.  For example, if $X=\R$ and $\mathcal{C}=\{(-\infty,t]:t\in\R\}$, then $VC(\mathcal{C})=1$; if $\mathcal{C}=\{[a,b]: a<b\in\R\}$, then $VC(\mathcal{C})=2$.  $VC$ dimension measures, in an intuitive sense, the complexity of a class of subsets.

The following lemma is a fundamental result in $VC$ theory, and we will use it several times.  A proof of the lemma and other details on the subject can be found in \cite{VC}. 

\begin{lm}[Sauer's Lemma]\label{T:Sauer}
 Let $\mathcal{C}$ be a class of subsets with $VC(\mathcal{C})=d<\infty$.  Then for any $k\geq d$, \[m^\mcc(k) \leq \left(\tfrac{ek}{d}\right)^d.\]
\end{lm}

For a class of functions $\mathcal{F}\subset\{f:X\to\{0,1\}\}$, denote by $\mcc_\mathcal{F}$ the set of subgraphs of functions in $\mathcal{F}$: $\mcc_\mathcal{F}=\{\{(x,t):t\leq f(x)\}:f\in\mathcal{F}\}.$  The $VC$ dimension of $\mathcal{F}$ is defined as $VC(\mcc_\mathcal{F})$, where this last quantity is the $VC$-dimension of a class of subsets of $X\times\R$.  It is worth noting that if $\mathcal{F}$ is the set of indicators of subsets in the class $\mcc$, $\mathcal{F}=\{1_C:C\in\mcc\}$, then $VC(\mathcal{F})=VC(\mcc).$  

It is well known in learning theory that empirical processes in the form of \eqref{E:emp} can be bounded via the $VC$-dimension of the indexing class.  Such results are often eponymously referred to as the ``$VC$ inequality'' after Vapnik and Chervonenkis, the pioneers of the theory.  In Section~\ref{S:noiseless} we use a version of the $VC$ inequality from \cite{panchenko}, which extends the $VC$ inequality to a more general case, when a class satisfies \emph{uniform entropy bounds}.  For a function $f$ and a probability $\Prb$, denote by $\Prb f$ the expectation $\int f\, d\Prb$.  For a class of binary functions $\mathcal{F}$ and a probability $\Prb$, the packing number $D(\mathcal{F},t,\Prb)$ is the cardinality of the largest subset $\mathcal{F}' \subset \mathcal{F}$ such that $\Prb|f-g|>t^2$ for all $f\not= g \in \mathcal{F}'$.  Finally, set \[D(\mathcal{F},t) := \sup_{\Prb} D(\mathcal{F},t,\Prb),\] where the supremum is taken over all discrete probabilities.  Then \cite[corollary 1]{panchenko} reads

\begin{thm}\label{T:panchenko}
Suppose \[\int_0^\infty \sqrt{\log D(\mathcal{F},t)}dt<\infty.\] Then there exists an absolute constant $K>0$ such that for any $u>0$ with probability at least $1-2e^{-u}$ for all $f\in \mathcal{F}$: 
\[\sum\limits_{k=1}^m \left(\Prb f -f(x_k)\right) \leq K \sqrt{m}\left( \sqrt{u\Prb f} + \int_0^{\sqrt{\Prb f}} \sqrt{\log(D(\mathcal{F},t))}  dt  \right).\]
\end{thm}

\subsection{Main VC Estimate}\label{S:VC}
This section is dedicated to the proof of Theorem~\ref{T:mainthmvc}.  We begin by computing the $VC$-dimension of all hemispheres, the case when $s=n$.

\begin{lm}\label{T:vchss} $VC(H_s^s)=s$.
 
\begin{proof}
 We first observe $H_s^s$ shattering the standard basis vectors $B=\{e_1,\dots,e_s\}$, and hence $VC(H_s^s)\geq s$.  Let $S\subset [s]$ and $B(S)=\{e_j:j\in S\}$.  Define $p=(p_1,\dots,p_s)$ by 
setting $p_j=1_S(j)-1_{S^c}(j)$.  Then $B(S)=B\cap H_p$.\\
\indent On the other hand, let $X=\{x_1,\dots,x_{s+1}\}$ be an arbitrary $(s+1)$-subset of $\mathbb{S}^{^{s-1}}$.  Without loss of generality, assume
\[x_{s+1} =\sum\limits_{k=1}^s \alpha_k x_k.\]
Set $A:=\{x_k:\alpha_k<0\}\cup \{x_{s+1}\}$; we'll see that for all $p\in\mathbb{R}^s$, $A \not= X\cap H_p$.  For any  $p$ such that $\langle p,x_k \rangle>0$ if $\alpha_k<0$ and $\langle p,x_k \rangle\leq 0$ if $\alpha_k\geq0$,
\begin{align*}
 \langle p, x_{s+1} \rangle &= \sum\limits_{k=1}^s \alpha_k \langle p,x_k \rangle\\
&= \sum\limits_{k:\alpha_k<0} \alpha_k \langle p,x_k \rangle + \sum\limits_{k:\alpha_k\geq 0} \alpha_k \langle p,x_k \rangle\\
&\leq 0.
\end{align*}
Therefore $H_s^s$ doesn't shatter $X$, so $VC(H_s^s)<s+1$.
\end{proof}
\end{lm}

We are now ready to estimate $VC(H^{n,s})$.  Let $d=VC(H^{n,s})\leq n$ and choose a subset $X=\{x^1,\dots,x^d\}$ of $\Sp^{n-1}$ shattered by $H^{n,s}$.  Fix an index set $S\in {[n] \choose s}$; for $x\in \Sp^{n-1}$ let $x_S=\sum_{j\in S} \avg{x,e_j}e_j$.  For any $B\subset \Sp^{n-1}$, let $B_S=\{b_S/\|b_S\|: b\in B \text{ and } b_S\not= 0\}$.  Notice that $|X_S|\leq d$, so by Lemmas~\ref{T:vchss} and \ref{T:Sauer}, \[\left|\{X_S\cap H_p:p\in \mathbb{S}^{n-1}_S\}\right| \leq \left(\tfrac{ed}{s}\right)^s.\]

The natural map $\{X\cap H_p:p \in \mathbb{S}^{n-1}_S\} \to \{X_S\cap H_p:p \in \mathbb{S}^{n-1}_S\}$ via $A\mapsto A_S$ is well-defined and surjective since $\sgn(\langle x,p\rangle)=\sgn(\langle x_S,p\rangle)$ for all $x\in X$ and $p\in \mathbb{S}^{n-1}_S$.  This map is also injective.  Suppose $A=X\cap H_p$ and $B=X\cap h_{p'}$ are distinct, for instance $a\in A\setminus B$ (hence $a_S \not= 0$).  If $a_S/\|a_S\| \in B_S$, then there is $b\in B$ such that $a_S/\|a_S\|=b_S/\|b_S\|$.  But then $\sgn(\langle a,p'\rangle) = \sgn(\langle b,p' \rangle),$ a contradiction.

It follows that \[\left|\{X\cap H_p: p\in \mathbb{S}^{n-1}_S\}\right|\leq \left(\tfrac{ed}{s}\right)^s,
\] and by the union bound, $2^d \leq {n \choose s}\left(\frac{ed}{s}\right)^s.$  After applying a familiar version of Stirlings approximation, ${n \choose s}\leq \left(\frac{en}{s}\right)^s,$ and some algebraic manipulation, we arrive at the inequality: \begin{equation}\label{preW}
-\log(2)\tfrac{d}{s} e^{-\log(2)\frac{d}{s}} \leq -\tfrac{\log(2)s}{e^2n}.
\end{equation}

To simplify further, we use the lower branch of the \textit{Lambert W function}, which is defined on $(\frac{-1}{e},0)$ by the relation $W_{-1}(x)e^{^{W_{-1}(x)}}=x$.  That is, $W_{-1}$ is the inverse of the map $x\mapsto xe^{x}$ restricted to $(-\infty,-1).$  We use the following lower bound of $W_{-1}$ to simplify \eqref{preW}.

\begin{lm}\label{T:lambo}
For all $-1/e<x<0,$ $W_{-1}(x)\geq \log(x^2)$.
\begin{proof}
Notice that $W_{-1}$ is decreasing, as is its inverse $W_{-1}^{-1}(x)=xe^x$.  Applying $W_{-1}^{-1}$ to each side of the equation in the statement and dividing by $x$, we find the equivalent: $x\log(x^2)\leq 1\,\text{ for all }\,\tfrac{-1}{e}<x<0.$  This holds since $x\mapsto x\log(x^2)$ is decreasing on $(-1/e,0)$ and $(\frac{-1}{e})\log(\frac{1}{e^2})=\frac{2}{e}<1.$
\end{proof}
\end{lm}

Applying the decreasing $W_{-1}$ to both sides of \eqref{preW} and using Lemma~\ref{T:lambo} gives: \[d\leq \tfrac{2}{\log(2)}s\log\left(\tfrac{ne^2}{s\log(2)}\right).\]


\section{The RIP of 1-Bit Embeddings}\label{S:noiseless}

This section proves Theorem~\ref{T:mainnoiseless2}.  Let $A\in \mathbb{R}^{m\times n}$ with rows $\{g_k\}_{k=1}^m$ drawn independently from the standard Gaussian distribution $\mathcal{N}(\textbf{0},I_n)$.  The Hamming distance between the images of two signals $x$ and $y$ under the 1-Bit embedding $\Phi_A$ is 
\begin{align*}
d_{\mch^m}(\Phi_Ax,\Phi_Ay) =  \frac{1}{m}\sum\limits_{k=1}^{m} \frac{1-\sgn\avg{x,g_k}\sgn\avg{y,g_k}}{2}.
\end{align*}
For $x,y\in\Sp^{n-1}$ we call $W_{x,y}:= H_x\triangle H_y $ (the symmetric difference of the two hemispheres) the \textit{wedge} associated to $x$ and $y$.  Notice $\sgn\avg{x,g_k} \not= \sgn\avg{y,g_k}$ if and only if $g_k$ is in the wedge $W_{x,y}$.  The Hamming distance above can be reformulated as \[d_{\mch^m}(\Phi_Ax,\Phi_Ay) = \frac{1}{m}\sum\limits_{k=1}^{m} 1_{W_{x,y}}(g_k).\]  The empirical processes framework suggests the sphere should be endowed with the distance $(x,y)\mapsto \prb(W_{x,y}).$  Fix $x,y\in \Sp^{n-1},$ let $g \sim \mathcal{N}(\textbf{0},I_n)$, and let $Z=(\avg{x,g}, \avg{y,g})^\top$; then $Z \sim \mathcal{N}\left(\textbf{0}, \Sigma\right)$ with \[\Sigma = \left[ \begin{array}{cc}
1 & \avg{x,y} \\
\avg{x,y} & 1 \end{array} \right].\]
It is a basic computation to find 
\begin{align}
\prb(W_{x,y}) &= \tfrac{1}{\pi\sqrt{1-\avg{x,y}^2}}\int_0^\infty\int_0^\infty \text{Exp}\left( \tfrac{2uv\avg{x,y}-u^2-v^2}{2-2\avg{x,y}^2} \right) du\,dv \label{E:prbd}\\[0.2em]
\nonumber &=\tfrac{1}{\pi} \arccos\big(\avg{x,y}\big).
 \end{align}
 This last quantity is the geodesic distance on the sphere that we denote by $d(x,y)$.  This brings our attention to the following object:
 \begin{align}
 \sup\limits_{x,y\in \Sp_s^{n-1}}\left| \frac{1}{m} \sum\limits_{k=1}^m 1_{W_{x,y}}(g_k) - d(x,y)\right|. \label{E:emp}
 \end{align}
The above formulation is paraphrased from \cite{BL}; this is the point at which our argument deviates.  To utilize the $VC$ theory for hemispheres developed in the previous section, we bound the $VC$-dimension of the class of ``sparse wedges'' $\mathcal{W}^{n,s}:=\{W_{x,y}:x,y\in \Sp_s^{n-1}\}.$

\begin{lm}
Let $\mcc$ be a class of subsets of $X$ with $VC(\mcc)=d<\infty.$  Let $\mcc\triangle\mcc =\{C\triangle C': C,C' \in \mcc\}.$  Then $VC(\mcc\triangle\mcc)\leq 10d.$
\begin{proof}
Let $B\subset X$ of size $m:=|B|$ to be prescribed later.  For a fixed pair $C,C'\in\mcc$, Notice that \[B\cap \left( C\triangle C'\right)= \left[(B\cap C)\setminus (B\cap C')\right] \cup \left[(B\cap C')\setminus (B\cap C)\right].\]  That is, $B\cap \left( C\triangle C'\right)$ is determined by $B\cap C$ and $B\cap C'$.  By Lemma~\ref{T:Sauer}, there are no more than $\left(\frac{em}{d}\right)^{2d}$ such pairs.  Taking $m\geq 10d$ yields $\left(\frac{em}{d}\right)^{2d} < 2^m.$
\end{proof}
\end{lm}
Along with Theorem~\ref{T:mainthmvc}, this lemma implies $VC(\mathcal{W}^{n,s})\lesssim s\log(n/s)$.  We use this $VC$-dimension estimate to bound the \textit{packing numbers} of the sparse wedges, $D(\mathcal{W}^{n,s},\eps,\Prb)$, which is the largest $d$ so that there exists $w_1,\dots, w_d\in \mathcal{W}^{n,s}$ with $\prb(w_i\triangle w_j)>\eps^2$ for all $i\not= j.$  General results bounding packing numbers via $VC$-dimension are well-known and the argument is standard; we include a proof in the current context for completeness.

\begin{prop}\label{T:packing}
For $0<\eps<1$, \[D(\mathcal{W}^{n,s},\eps,\Prb) \lesssim \left(\tfrac{1}{\eps^2}\right)^{VC(\mathcal{W}^{n,s})+1}.\]
\begin{proof}
Fix $0<\eps<1$.  Let $d=D(\mathcal{W}^{n,s},\eps,\Prb)$ and let $w_1 \dots, w_d$ such that $\prb(w_i\triangle w_j)>\eps^2$ for all $i\not= j$.  Let $\{X_k\}_{k=1}^n$ be independent and identically distributed on the sphere with law $\prb$, where $n$ will be determined later.  Notice that $w_i \cap \{X_k\} \not= w_j \cap \{X_k\}$ if and only if $(w_i\triangle w_j) \cap \{X_k\}$ is nonempty.  Thus the probability that there is $i\not=j$ such that $w_i\cap\{X_k\}=w_j\cap\{X_k\}$ is no more than 
\begin{align*}
\tbinom{d}{2}\,{\raisebox{0.4 em}{$\max\limits_{1\leq i\not= j\leq d}$}}\,\,\prb\left(w_i\cap\{X_k\}=w_j\cap\{X_k\}\right) &= \tbinom{d}{2}\,{\raisebox{0.4 em}{$\max\limits_{1\leq i\not= j\leq d}$}}\,\,\left(1-\prb(w_i \triangle w_j)\right)^n\\
&< \tbinom{d}{2} (1-\eps^2)^n\\
&< d^2 e^{-n\eps^2}\\
&= e^{2\log(d) - n\eps^2}.
\end{align*}
Now we take $n=\ceil*{\frac{2\log(d)+1}{\eps^2}}$ so the above probability is less than one, hence there is a deterministic $X=\{x_k\}_{k=1}^n$ so that $d=|\{w_j\cap X:j\in [d]\}|.$  Let $v=VC(\mathcal{W}^{n,s}).$  Employing Lemma~\ref{T:Sauer}, there is $K_v>0$ such that 
\begin{align*}
d &\leq K_v \left( \tfrac{2\log(d)+2}{\eps^2}\right)^v.\\
\end{align*}
Choose $d_0$ large enough so that for $d>d_0,$ $(2\log(d)+2)^{v+1}<d^{1/v}$.  This yields \[d\leq \max\{d_0,K_v^{\frac{v+1}{v}}\} \left(\tfrac{1}{\eps^2}\right)^{(v+1)}.\]
\end{proof}
\end{prop}

Notice that the bound in Proposition~\ref{T:packing} holds uniformly over all probabilities on the sphere.  This fact allows us to use a version of the entropy integral in the final stage of our argument.  Recall Theorem~\ref{T:panchenko}.  Adapted to our current setting, we have the following corollary:

\begin{cor}\label{T:pancor1}
There exists an absolute constant $K>0$ such that for any $u>0$ with probability at least $1-2e^{-u}$ for all $W_{x,y}\in \mathcal{W}^{n,s}$: 
\[\sum\limits_{k=1}^m \left(d(x,y)-1_{W_{x,y}}(g_k)\right) \leq K \sqrt{m}\left( \sqrt{ud(x,y)} + \int_0^{\sqrt{d(x,y)}} \sqrt{\log(D(\mathcal{W}^{n,s},t))}  dt  \right).\]
\end{cor}

We adjust this result in two ways to produce the main results of this section.  First, increase the right side of the inequality by replacing all distances with one.  Now that the bound is uniform over pairs of signals in $\Sp_s^{n-1}$, we observe \[\sup_{x,y\in\Sp_s^{n-1}} \sum\limits_{k=1}^m \left(d(x,y)-1_{W_{x,y}}(g_k)\right) = \sup_{x,y\in\Sp_s^{n-1}} \left|\sum\limits_{k=1}^m 1_{W_{x,y}}(g_k)-d(x,y)\right|.\]  This is because $1_{W_{-x,y}} = 1-1_{W_{x,y}}$ (a.s.), and $d(-x,y)=1-d(x,y)$.  Thus we have:

\begin{cor}
There exists an absolute constant $K>0$ such that for any $u>0$ with probability at least $1-2e^{-u}$, 
\[\sup_{x,y\in\Sp_s^{n-1}} \frac{1}{m}\left|\sum\limits_{k=1}^m 1_{W_{x,y}}(g_k)-d(x,y)\right| \leq \frac{K}{\sqrt{m}}\left( \sqrt{u} + \int_0^1 \sqrt{\log(D(\mathcal{W}^{n,s},t))}  dt  \right).\]
\end{cor}

After applying the uniform entropy bounds of Proposition~\ref{T:packing} in the above corollary and setting $u=\log(2/\eps)$, Theorem~\ref{T:mainnoiseless2} is immediate.


\section{The RIP of Noisy 1-Bit Embeddings}\label{S:noisy}

\subsection{Noisy RIP with the distorted metric on the sphere}
This section proves Theorem~\ref{T:mainnoise}.  We again consider $A\in \mathbb{R}^{m\times n}$ with rows $\{g_k\}_{k=1}^m$ drawn independently from the standard Gaussian distribution $\mathcal{N}(\textbf{0},I_n)$.  We are now interested in the case of additive white noise prior to quantization; let $\eta \sim \mathcal{N}(\textbf{0},\sigma^2 I_m).$  Then the Hamming distance between the images of two signals under the 1-Bit embedding $\Phi_A^\eta$ is
\begin{align}\label{E:ham}
d_{\mch^m}(\Phi^\eta_Ax,\Phi^\eta_Ay) =  \frac{1}{m}\sum\limits_{k=1}^{m} \frac{1-\sgn(\avg{x,g_k} +\eta_k)\sgn(\avg{y,g_k}+\eta_k)}{2}.
\end{align}

Fix $x,y\in \Sp^{n-1}.$  Let $g \sim \mathcal{N}(\textbf{0},I_n)$ and $\mu \sim \mathcal{N}(0,\sigma^2)$ be independent.  Then ${\avg{x,g}+\mu \choose \avg{y,g}+\mu}$ is a Gaussian vector with covariance matrix \[\left[ \begin{array}{cc}
1+\sigma^2 & \avg{x,y}+\sigma^2 \\
\avg{x,y}+\sigma^2 & 1+\sigma^2 \end{array} \right].\]
A computation similar to \ref{E:prbd} yields
\begin{align*}
\prb \big( \sgn(\avg{x,g} +\mu)\sgn(\avg{y,g}+\mu)=-1 \big) &=\tfrac{1}{\pi} \arccos\left(\tfrac{\avg{x,y}+\sigma^2}{1+\sigma^2}\right).
 \end{align*}
This last quantity is $d^\sigma(x,y)$, defined in \eqref{e:xdsigma}.  We'll see soon that $d^\sigma$ is in fact a metric; this is the distance with which $\Sp_s^{n-1}$ is naturally endowed in the presence of additive white noise.  The object in the $(s,\delta,
\sigma)$-RIP that we aim to bound is \[\sup\limits_{x,y\in \Sp_s^{n-1}}\Big| d_{\mch^m}(\Phi^\eta_Ax,\Phi^\eta_Ay) - d^\sigma(x,y)\Big|.\]

\begin{center}
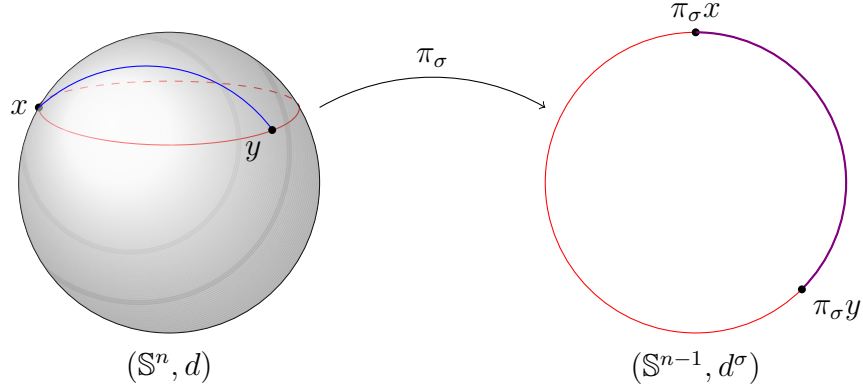
\begin{figure}
\centering
\begin{tikzpicture}
    \draw[red] (-1.73205,1) arc (180:360:1.73205cm and 0.5cm);
    \draw[red, dashed] (-1.73205,1) arc (180:0:1.73205cm and 0.35cm);
    \draw[->] (2,1) arc (120:60:3cm);
    \draw[red] (7,0) circle (2cm);
    \draw (7,-2.1) node[anchor=north] {$(\Sp^{n-1},d^\sigma)$};
    \fill (-1.73205,1) circle[radius=1.5pt];
    \draw (-1.73205,1) node[anchor=east] {$x$};
    \fill (7,2) circle[radius=1.5pt];
    \draw (7,2) node[anchor=south] {$\pi_\sigma x$};
    \draw (0,0) circle (2cm);
    \draw (0,-2.1) node[anchor=north] {$(\Sp^n,d)$};
    \draw (3.5,1.4) node[anchor=south] {$\pi_\sigma$};
    \shade[ball color=gray!10!white,opacity=0.20] (0,0) circle (2cm);
    
    \draw[blue] (-1.73205,1) arc (132.5:36:2.1cm);
	\fill (1.37,.7) circle[radius=1.5pt];
	\draw (1.37,.7) node[anchor=north east] {$y$};
	\draw[violet, thick] (7,2) arc (90:-45:2cm);
	\fill (8.412,-1.42) circle[radius=1.5pt];
	\draw (8.412,-1.42) node[anchor=north west] {$\pi_\sigma y$};
	
\end{tikzpicture}
\caption{If $\pi_\sigma:\{p\in\Sp^n:\avg{p,e_{n+1}}=\sigma\} \to \Sp^{n-1}$ is the normalization of the projection onto the first $n$ coordinates, then $d^\sigma(\pi_\sigma x,\pi_\sigma y)=d(x,y).$} \label{fig:dsigma}
\end{figure}
\end{center}

Appealing to the methods in Section~\ref{S:noiseless}, we rewrite the additive noise as an inner product by increasing the Gaussian measurements' dimension by one and lifting the sphere to one higher dimension by padding with $\sigma^2$.  Introduce the following notation:\begin{align*}
x_\sigma &=\tfrac{1}{\sqrt{1+\sigma^2}}(x^1,\dots,x^n,\sigma)\in \Sp^n_{s+1}.
\end{align*}
Let $h=(g^1,\dots,g^n,\frac{1}{\sigma}\mu)$ and notice $\avg{x_\sigma,h} =\tfrac{1}{\sqrt{1+\sigma^2}}(\avg{x,g}+\mu)$ and $h\sim \mathcal{N}(\textbf{0},I_{n+1}).$  Denote by $W^\sigma_{x,y}$ the wedge in $\Sp^{n}$ relative to $x_\sigma$ and $y_\sigma$, i.e., \[W^\sigma_{x,y}:=H_{x_\sigma}\triangle H_{y_\sigma}.\]  Then $\sgn(\avg{x,g_k}+\eta_k) \not= \sgn(\avg{y,g_k}+\eta_k)$ if and only if $h_k:=(g^1_k,\dots,g^n_k,\frac{1}{\sigma}\eta_k) \in W^\sigma_{x,y}$.  The Hamming distance in \eqref{E:ham} can be reformulated as \[d_{\mch^m}(\Phi^\eta_Ax,\Phi^\eta_Ay) = \frac{1}{m}\sum\limits_{k=1}^{m} 1_{W^\sigma_{x,y}}(h_k).\]

Furthermore, notice that \[\left\langle x_\sigma,y_\sigma\right\rangle=\tfrac{\avg{x,y}+\sigma^2}{1+\sigma^2},\] hence $d^\sigma(x,y)=d\left( x_\sigma,y_\sigma\right),$ where we abuse notation to allow $d\left( \cdot,\cdot\right)$ to denote the normalized geodesic distance on $\Sp^{n}$; see Figure~\ref{fig:dsigma} for an illustration.  It is now apparent that $d^\sigma$ is indeed a metric on $\Sp^{n-1}$.

This brings our attention to the following object:
 \[\sup\limits_{x,y\in \Sp_s^{n-1}}\left| \frac{1}{m} \sum\limits_{k=1}^m 1_{W^\sigma_{x,y}}(h_k) - d^\sigma\left(x,y\right)\right|,\]
 where the row vectors $h_k$ are independently drawn from $\mathcal{N}(\textbf{0},I_{n+1}).$  At this point, the argument of Section~\ref{S:noiseless} applies so long as we can estimate the $VC$-dimension of \[\mathcal{W}^{n,s}_\sigma:=\{W_{x,y} \in \mathcal{W}^{n+1,s+1}: x^{n+1}=\sigma=y^{n+1}\}.\]  Since $\mathcal{W}^{n,s}_\sigma \subset \mathcal{W}^{n+1,s+1}$, it is clear that
 \[VC(\mathcal{W}_\sigma^{n,s})\leq VC(W^{n+1,s+1}) \lesssim (s+1)\log\left(\tfrac{n+1}{s+1}\right) \lesssim s\log(n/s).\]
 
\subsection{Noisy RIP with geodesic metric on the sphere}\label{S:quieting}
The deviation of $d^\sigma$ from the geodesic distance is exaggerated at antipodes.  That is, for any $x$ and $y$ on the sphere, $|d(x,y)-d^\sigma(x,y)|\leq d(x,-x)-d^\sigma(x,-x)$.  In what is to come, all suprema are over $x,y\in\mathbb{S}^{n-1}_s.$  If one prefers a bound of the form \[\sup\left| d_{\mathcal{H}^m}(\Phi^\eta_Ax,\Phi^\eta_Ay)-d(x,y)\right| \leq \delta,\] it is enough for 
\begin{align*}\sup\left| d_{\mathcal{H}^m}(\Phi^\eta_Ax,\Phi^\eta_Ay)-d(x,y)\right| &\leq \sup\left| d_{\mathcal{H}^m}(\Phi^\eta_Ax,\Phi^\eta_Ay)-d^\sigma(x,y)\right| + \sup\left| d^\sigma(x,y)-d(x,y)\right|\\
&\leq \sup\left| d_{\mathcal{H}^m}(\Phi^\eta_Ax,\Phi^\eta_Ay)-d^\sigma(x,y)\right| + 1-\tfrac{1}{\pi}\arccos\left(\tfrac{\sigma^2-1}{\sigma^2+1}\right)\\
&\leq \delta.
\end{align*}
Corollary~\ref{T:quieting} follows easily from Theorem~\ref{T:mainnoise} and this observation.

\subsection*{Acknowledgements}The author would like to thank Michael Lacey for suggesting the project and pointing out the utility of the $VC$ theory, as well as for useful conversation on the subject and layout of the article.


\end{document}